\documentclass[12pt]{article}
\usepackage{amsthm,amstext,amscd,latexsym,amsmath,amsfonts,amssymb,graphicx}
\def\dj{d\kern-.30em\raise1.25ex\vbox{\hrule width .3em height
.03em}}
\def\Dj{D\kern-.70em\raise0.75ex\vbox{\hrule width .3em height
.03em}\kern.10em}
\newcommand{\dwedge}{\,\dot\wedge\,}

\newcommand{\openone}{\text{1\kern-0.25em I}}

\newcommand{\Z}{\mathbb{Z}}

\newcommand{\R}{\mathbb{R}}\newcommand{\N}{\mathbb{N}}

\newcommand{\Cl}{\mathcal{C}\kern -0.15em\ell}
\newcommand{\ie}{\textit{i.e.}\,}
\newcommand{\alg}{\text{alg}}\newcommand{\cog}{\text{cog}}
\newcommand{\tr}{\text{tr}}

\newcommand{\End}{\text{End\,}}
\newcommand{\Hom}{\text{Hom\,}}
\newcommand{\id}{\text{id}}
\newcommand{\im}{\text{im}}
\newcommand{\half}{\textstyle{\frac{1}{2}}}

\theoremstyle{definition}
\newtheorem{dfn}{Definition}[section]
\newtheorem{exm}[dfn]{Example}
\newtheorem{prob}[dfn]{Problem}
\newtheorem{Lemma}[dfn]{Lemma}
\newtheorem{Mnthm}[dfn]{Main Theorem}
\newtheorem{thm}[dfn]{Theorem}
\newtheorem{rem}[dfn]{Remark}

\begin{document}\title{\vspace{-2truecm}CLIFFORD HOPF GEBRA FOR\\
TWO-DIMENSIONAL SPACE\thanks{Submitted to `Miscellanea Algebraica', 
Waldemar Korczy\'nski and Adam Obtu{\l}owicz, Editors, Akademia 
\'Swi{\c e}tokrzyska, Kielce, Poland.}}
\author{Bertfried Fauser\\Universit\"at Konstanz, 
Fachbereich Physik, Fach M678\\ 
D-78457 Konstanz\\Bertfried.Fauser@uni-konstanz.de 
\and Zbigniew Oziewicz\thanks{Supported by el Consejo Nacional de
Ciencia y Tecnolog\'{\i}a (CONACyT) de M\'exico, grant \# 27670 E
(1999-2000), and by UNAM, DGAPA, Programa de Apoyo a Proyectos de
Investigaci\'on e Innovaci\'on Tecnol\'ogica, Proyecto IN-109599
(1999-2002). This work was also supported partially by Universit\"at
Konstanz (Germany). Zbigniew Oziewicz is a member of Sistema Nacional de
Investigadores, M\'exico, No. de expediente 15337.}\\
Universidad Nacional Aut\'onoma de M\'exico\\
Facultad de Estudios Superiores Cuautitl\'an\\
Apartado Postal \# 25, C.P. 54700 Cuautitl\'an Izcalli,\\ 
Estado de M\'exico\\oziewicz@servidor.unam.mx}
\date{Submitted November 22, 2000}\maketitle
\vspace{-1truecm}
\begin{abstract}
A Clifford algebra $\Cl(V,\eta\in V^*\otimes V^*)$
jointly with a Clifford cogebra $\Cl(V,\xi\in V\otimes V)$ is said to
be a Clifford biconvolution $\Cl(\eta,\xi).$ We show that a Clifford
biconvolution for $\dim_\R\Cl=4$ does possess an antipode iff
$\det(\id-\xi\circ\eta)\neq 0.$ An antipodal Clifford biconvolution is
said to be a Clifford Hopf gebra.\\[1ex]
\noindent\textbf{2000 Mathematics Subject Classification:} 
15A66 Clifford algebra, 16W30 Coalgebra, bialgebra, Hopf algebra\\
\noindent\textbf{2000 PACS:} 02.10.Tq Associative rings and algebras\\
\noindent\textbf{Keywords:} cliffordization, Clifford algebra, Clifford
cogebra, antipode, Hopf gebra, Clifford bigebra, Gra{\ss}mann algebra
\end{abstract}\tableofcontents

\section{Introduction} After Bourbaki [1989 \S 11] we use
\textit{cogebra, bigebra} and Hopf \textit{gebra} instead of coalgebra,
bialgebra and Hopf algebra.

Let $C$ be an $\R$-space and $C^*$ be an $\R$-dual $\R$-space.
If $C$ is an $\R$-cogebra and $A$ is an $\R$-algebra then the $\R$-space
$A\otimes C^*$ inherits the structure of an $\R$-algebra with a
convolution product: this is a convolution $\R$-algebra, an
$\R$-convolution for short. A dual $\R$-space $C\otimes A^*$ inherits a
structure of an $\R$-cogebra with coconvolution coproduct: this is a
coconvolutional $\R$-cogebra, an $\R$-coconvolution for short.

In particular if an $\R$-space $V$ carries an $\R$-biconvolution algebra
\& cogebra structure, then do also the $\R$-spaces $\End V\simeq V\otimes
V^*,\;\End V^*,$ as well as all iterated $\R$-spaces $\End\ldots\End V$
inherit also $\R$-biconvolution algebra \& cogebra structures.

If the $\R$-space $C$ is an $\R$-cogebra with a coproduct
$\triangle:C\rightarrow C\otimes C,$ then $C^*$
is an $\R$-algebra with product $\triangle^*:C^*\otimes C^*\rightarrow
C^*.$ If an $\R$-space $A$ is a finite dimensional $\R$-algebra having a
binary product $m:A\otimes A\rightarrow A,$ then an $\R$-dual $\R$-space
$A^*$ (or $\Z$-graded dual in the case $A$ is not a finite dimensional
$\R$-space) is an $\R$-cogebra with a binary coproduct
$m^*:A^*\rightarrow A^*\otimes A^*.$

However there are several important situations (free tensor algebra,
exterior algebra, Clifford algebra, Weyl algebra, ...) where the dual
space of an algebra is also an algebra in a natural way by construction.
If this is the case, then by the above ($\Z$-graded) duality, both
mutually dual $\R$-spaces carry both structures, algebra \& cogebra, and
therefore we have a dual pair of $\R$-biconvolutions.

An unital and associative convolution possessing an (unique) antipode is
said to be a Hopf gebra or an \textit{antipodal convolution}
(Definition \ref{dfn2.2}). This terminology has been introduced
in [Oziewicz 1997, 2001; Cruz \& Oziewicz 2000] and is different
from Sweedler's [1969 p. 71]. A general theory of the
finite-dimensional antipodal and antipode-less biconvolutions
(convolutions and coconvolutions) has been initiated in [Cruz \&
Oziewicz 2000]. 

Nill in 1994 and B\"ohm \& Szlach\'anyi in 1996 introduced weak bigebras
and weak Hopf gebras with antipode $S$ defined as the Galois
connection with respect to the binary convolution $*$ which does not
necessarily needs to be unital [Nill 1998, Nill et all. 1998],
$$
\id *S*\id=\id,\quad S*\id*S=S.
$$ 

If $\eta\in V^{*\otimes 2}$ is invertible then $\eta^{-1}\in V^{\otimes
2}.$ If $\Cl(V^*,\xi\in V^{\otimes 2})$ is a Clifford $\R$-algebra, then
a dual $\R$-space $\Cl(V,\xi)\equiv\{\Cl(V^*,\xi)\}^*$ is a Clifford
$\R$-cogebra. It was shown in [Oziewicz 1997, and ff.] that a Clifford
convolution $\Cl(\eta,\eta^{-1})\equiv\Cl(V,\eta\in V^{*\otimes
2},\eta^{-1}\in V^{\otimes 2})$ is antipode-less.

The aim of this paper is to show that a Clifford convolution
$\Cl(\eta,\xi)\equiv\Cl(V,\eta,\xi)$ for $\dim_\R V=2,$
\ie for $\dim_\R(V^\wedge)=\dim_\R\Cl(V)=4,$ does posses an antipode
iff $\det(\id-\xi\circ\eta)\neq 0$ (Main Theorem 4.1).
Applications to physics have been proposed in [Fauser 2000b].

\section{Biconvolution and antipode}
\begin{figure}[ht]\hfil 
\includegraphics[height=2truecm]{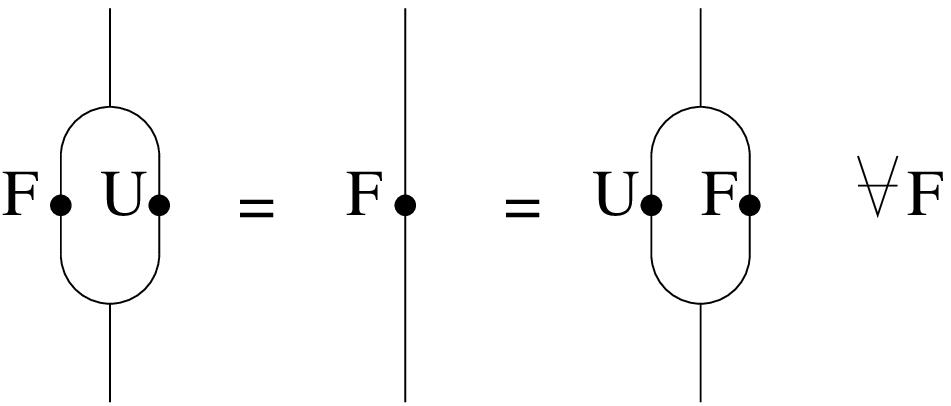}\hfil 
\caption{Unital convolution. We assume
$U=u\circ\epsilon\in\End V.$}
\label{conv_unit}\end{figure}   
\begin{figure}[ht]\hfil\includegraphics[height=2truecm]{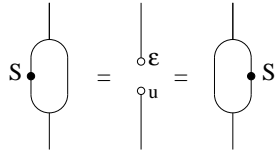}
\hfil \caption{Axioms for the antipode $S.$}\label{antipode}
\end{figure}
Let the convolution be unital, Fig. \ref{conv_unit}. This is the case if
an $\R$-algebra $A$ is unital with unit $u:\R\rightarrow A$ and an
$\R$-cogebra $C$ is counital with counit $\epsilon:C\rightarrow\R.$
\begin{dfn}[Sweedler 1969 p. 71, Zakrzewski 1990 p. 357]\label{dfn2.1}
The convolutive inverse, w.r.t. the convolutive unit 
$U=u\circ\epsilon \in \End\,V$, of the identity map on $V$, 
Fig. \ref{antipode}, is said to be an \textit{antipode},
$S\equiv(\id)^{-1}$.\end{dfn}
If an antipode exists w.r.t. an unital associative convolution it must
be unique.

An antipodal biconvolution defines a unique crossing as given in Fig.
\ref{crossing} in the last Section, see [Oziewicz 1997, 2001; Cruz
\& Oziewicz 2000]. Using the axioms of the antipode, Fig. \ref{antipode},
and biassociativity, this crossing is equivalent to the algebra 
homomorphism between algebra and crossed algebra, as well as, to the 
cogebra homomorphism from crossed cogebra to cogebra. A proof is given 
in [Cruz \& Oziewicz 2000].

\begin{exm}[Gra{\ss}mann Hopf gebra] The Gra{\ss}mann wedge product,
the Gra{\ss}mann coproduct and the unique antipode, as given by
Sweedler [1969, Ch. XII] and extended by Woronowicz [1989, \S 3], see 
Section 3.1 below for details, gives the Gra{\ss}mann Hopf gebra.
\end{exm} This motivates the following definition:
\begin{dfn}\label{dfn2.2} An unital and associative convolution
possessing a (unique) antipode is said to be a Hopf gebra or an
antipodal convolution.\end{dfn}

\begin{figure}[ht]\hfil
\includegraphics[height=2truecm]{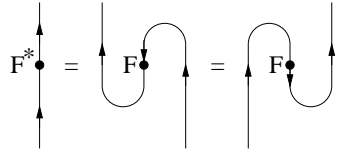}\hfil
\caption{Transposition [Lyubashenko 1995, fig 1 on p. 249].}
\label{dual_op}\end{figure}

\begin{figure}[ht]\hfil
\includegraphics[height=2truecm]{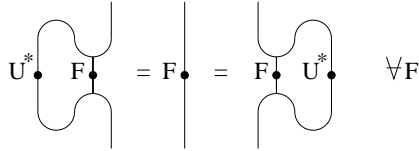}\hfil
\caption{Coconvolutive counit $U^*\in\End V^{*\wedge}.$}
\label{conv_counit}\end{figure}

A closed structure is given by the evaluation and coevaluation 
[Kelly \& Laplaza 1980; Lyubashenko 1995]:
$$\begin{CD}V^*\otimes V@>{\text{ev}}>>\R,\\
V^*\otimes V@<<{\text{coev}}<\R.\end{CD}$$
Following [Lyubashenko 1995, fig. 3, p. 250] the evaluation is
represented by cup and coevaluation by cap. The evaluation intertwines
the transposed endomorphisms $F^*\in\End V^*$ with $F\in\End V$,
Fig. \ref{dual_op}. Up and down arrows indicate the spaces $V^*$ and
$V$.
A coconvolution is counital if e.g. there exits a dual counit
$u^*:V^{*\wedge}\mapsto\R^*$ and a dual unit 
$\epsilon^*:\R^*\mapsto V^{*\wedge}$, Fig. \ref{conv_counit}.
\begin{figure}[ht]\hfil
\includegraphics[height=2truecm]{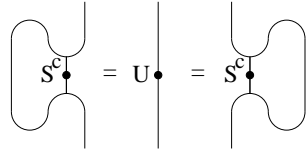}\hfil
\caption{Coconvolutive coantipode $S^c.$}
\label{conv_coantipode}\end{figure}
We have to use the dual coconvolution counit $U^*\in\End V^{*\wedge}$
in Fig. \ref{conv_counit}.

\section{Rota \& Stein's cliffordization} 
\subsection{Exterior Gra{\ss}mann Hopf gebra} The Gra{\ss}mann Hopf
gebra was constructed by Sweedler [1969, Ch. XII] factoring the
couniversal shuffle tensor Hopf gebra Sh$V$ by the \textit{switch},
$s(x\otimes y)\equiv y\otimes x$. This construction was generalized to
any braid by Woronowicz [1989, \S 3, p. 154]. An exterior Hopf gebra
can be defined in terms of an unique braid dependent homomorphism of
universal tensor Hopf gebra into couniversal tensor Hopf gebra and this
implies that an exterior Hopf gebra is couniversal and braided
[Oziewicz, Paal \& R\'o\.za\'nski 1995 \S 8; R\'o\.za\'nski 1996;
Oziewicz 1997 p. 1272-1273].
\begin{figure}[ht]\hfil
\includegraphics[height=2.5truecm]{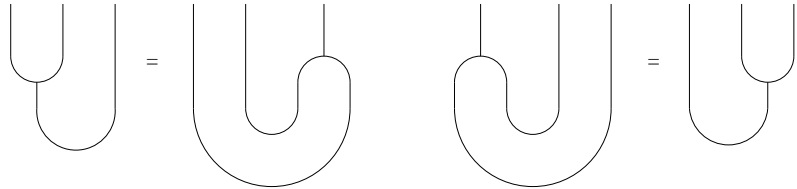}\hfil
\caption{Left-right product - coproduct duality. 
Cup's are either `ev' or $\eta^\wedge.$}
\label{duality}\end{figure}

\subsection{Cliffordization} 
\begin{figure}[ht]\hfil
\includegraphics[height=2truecm]{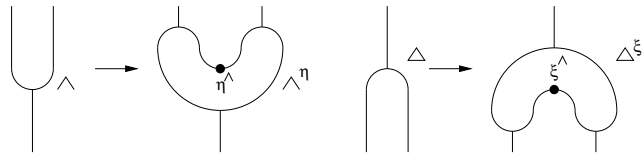}\hfil
\caption{Bicliffordization: the sausage graphs.}
\label{cliffordization}\end{figure}

The tensors $\eta,\eta^T\in V^*\otimes V^*$
are said to be scalar products on $V$ or coscalar products on $V^*$
($\eta^T$ is the transpose of $\eta$). The tensors $\xi,\xi^T\in
V\otimes V$ are said to be scalar products on $V^*$ and coscalar on $V$.
In particular $\half(\eta+\eta^T)$ is the symmetric part of
$\eta.$ The scalar products are displayed by decorated (or labelled)
cups and coscalar products by decorated caps, see Fig.
\ref{cliffordization}.
 
Rota and Stein [1994] introduced the Clifford product as a deformation 
of exterior biconvolution. This deformation process was called 
\textit{cliffordization}. A cliffordization introduces an internal
loop in a binary product having two inputs and one output, 
employing the $\eta^\wedge$-cup (or $\xi^\wedge$-cap for the coproduct), 
Fig. \ref{cliffordization}.

Clifford biconvolution was defined in [Oziewicz 2001] as the
$(\eta,\xi)$-bicliffordization of an exterior biconvolution. In
Sweedler's notation,
\begin{gather*}\wedge^\eta(x\otimes y):=x_{(1)}\wedge\eta^\wedge
(x_{(2)}\otimes y_{(1)})\wedge y_{(2)},\\
\Delta^\xi 1=\xi^\wedge\quad\text{and}\quad
\triangle^\xi x:=x_{(1)}\wedge\xi^\wedge\wedge x_{(2)}.\end{gather*}
If the product coproduct duality of Fig. \ref{duality} is used 
with cup as an evaluation, then every product on $V^\wedge$ induces a
coproduct on $V^{*\wedge}$ and vice versa. If $\eta^\wedge$-cup's and 
$\xi$-cap's are used, one gets a correlation between products on 
$V^\wedge$ and coproducts on $V^\wedge$. 

A Clifford $\R$-algebra together with a Clifford $\R$-cogebra on
$V^\wedge,$ $\Cl(\eta,\xi)\equiv\Cl(V,\wedge^\eta,\Delta^\xi),$ is said
to be a Clifford $\R$-convolution. It was shown in [Oziewicz 1997] that
a Clifford $\R$-convolution for $\xi\circ\eta=\id$ and for
$\eta\circ\xi=\id$ is antipode-less. An antipode-less Clifford
convolution $\Cl(\eta,{\eta^{-1}})$ for an {\it invertible} tensor
$\eta$ \textit{cannot} be a deformation of an exterior Gra{\ss}mann 
Hopf gebra.

\begin{dfn}The tensors $\eta$ and $\xi$ are said to be dependent if
$0 \neq A\in \End\,V$ and $0 \neq B\in \End\, V^*$ exist such that one
of the following relations hold, $$\xi \,=\, A \circ \eta^{-1}\circ B,
\qquad\eta \,=\, B \circ \xi^{-1} \circ A.$$\end{dfn}
If the tensors $\eta\,\&\,\xi$ are independent then the Clifford
product and coproduct are defined independently by Rota \& Stein's
\textit{deformation}. 

Cliffordization and cocliffordization of the Gra{\ss}mann convolution 
does not change the convolutive unit $U=u\circ\epsilon.$ However since 
the $\Z$-grading is changed due to the skewsymmetric parts of $\eta$
and $\xi,$ the counit is no longer the projection onto
$\R\subset V^{\dwedge}$ [Fauser 1998-2000a].

\section{The Clifford antipode} 
\begin{Mnthm}\label{MnThm} A Clifford biconvolution $\Cl(\eta,\xi)$ is
a Clifford Hopf gebra iff $\mu\equiv\det(\id-\xi\circ\eta)\neq 0.$ Then 
\begin{align*}
(i)&\quad\mu S\vert_{V}=-\id_V.\\
(ii)&\quad\det(\mu S)=(-1)^{\dim V}.\\ 
(iii)&\quad \tr(\mu S)=\half\tr\{(\eta-\eta^T)\circ(\xi-\xi^T)\}.\\
(iv)&\quad\text{The minimal polynomial of $\mu S$ is}\\
&\quad(\mu S+1)[(\mu S-1)^{\dim V}-\tr(\mu S)\cdot\mu S].\end{align*}
\end{Mnthm}
\begin{proof} Theorem 4.1 was proved for $\dim_\R V=1$ in [Oziewicz
1997]. If $\dim_\R V=n$ then for $A\in\End V,$
$$\det(\lambda\cdot\id-A)=\det(-A)+\ldots+\lambda^{n-2}\half
\left\{(\tr A)^2-\tr(A^2)\right\}-\lambda^{n-1}\tr A +\lambda^n.$$
Therefore $\det(\id-\eta\circ\xi)=\det(\id-\xi\circ\eta).$

The Clifford antipode $S\equiv S(\eta,\xi)\in\End V^\wedge$
is computed from Fig. \ref{antipode}. We give the proof for 
$\dim_\R V=2$ only. The general case will be treated elsewhere. 
Let $r,s,t,u,v,z\in\R$ and
\begin{gather}
\eta e_1=+r^2\epsilon^1+t\epsilon^2,\qquad\xi\epsilon^1=
+u^2e_1+ze_2,\nonumber\\
\eta e_2=-t\epsilon^1-s^2\epsilon^2,\qquad\xi\epsilon^2=-ze_1-v^2e_2.
\label{eta}
\end{gather} 
Then we find together with (i) of the main theorem:
\begin{align*}
\mu\cdot S\,1\quad&=1+4zt+2t\,e_1\wedge e_2,\\
\mu S(e_1\wedge e_2)&=\quad 2z\quad+\quad e_1\wedge e_2.
\end{align*}
An action of $g\in GL(V)$ on tensors
$\eta\in V^*\otimes V^*$ and $\xi\in V\otimes V$ is given as follows
$$
\eta\mapsto g^T\circ\eta\circ g,\qquad\xi\mapsto g\circ\xi\circ g^T.
$$
It would be desirable to study $GL(V)$-orbits on $(V^*\otimes V^*)
\times(V\otimes V)$ and present full classification of all orbits in
terms of invariants. However, this topic exceeds the scope of this
paper and will be presented elsewhere.
\end{proof}
\begin{rem}
Theorem \ref{MnThm} solves and improves Conjecture 2.2 posed
in [Oziewicz 1997, p. 1270].
\end{rem}
\begin{rem} The above Clifford Hopf gebra includes as a particular 
case for $\xi=0$ the construction made by {\Dj}\,ur{\dj}evich
[1994].
\end{rem}
\begin{rem}In the case $\det\eta\neq 0$ one can take $\xi=\eta^{-1}$.
Also if $\det\xi\neq 0$ one can choose $\eta=\xi^{-1}$. In order to
prevent such possibilities we need to supplement the definition of 
the Clifford Hopf gebra with the extra conditions:
\begin{gather*}\det(\id-\xi\circ\eta)\neq 0,\quad
\det\xi=\det\eta=0.
\end{gather*}
The antisymmetric tensor $F\equiv\half(\eta-\eta^T)$ can be adjusted in
such a way that $\det\eta=0$ with invertible symmetric tensor
$g\equiv\half(\eta+\eta^T).$
We found a Clifford antipode for $\det\eta=0=\det \xi,$
\begin{gather*}\det(\id_V-\xi\circ\eta)=-\tr(\xi\circ\eta)+1,\\
t=\pm rs,\quad z=\pm uv;\qquad\tr(\xi\circ\eta)=(ur\pm sv)^2\neq 1.
\end{gather*}\end{rem}
The following matrices for $rs\neq 0$ represent the same tensor from
$V^{*\otimes 2}$ or from $V^{\otimes 2}$ with respect to the different
bases, these matrices are on the same $GL(2,\R)$ orbit,
$$\begin{pmatrix}1&1\\-1&-1\end{pmatrix}\simeq
\begin{pmatrix}r^2&rs\\-rs&-s^2\end{pmatrix}\simeq
\begin{pmatrix}r^2&2r^2s\\0&0\end{pmatrix}.$$

An antipode for regular scalar and coscalar tensors can be found also.

\section{An antipode-less Clifford bigebra} 
A Clifford biconvolution $\Cl(\eta,\xi)$ is antipode-less if 
$\det(\id-\xi\circ\eta)=0$. In particular this is the case if 
$\xi=\eta^{-1}.$ It was shown [Oziewicz 1997] that 
$\wedge^\eta\circ\Delta^{\eta^{-1}}=(\dim\Cl)\cdot\id_{\Cl}$.
\begin{prob} 
What axioms for Clifford biconvolution implies the condition 
$\det(\id-\xi\circ\eta)=0$? In particular, does a braid
exists for which such a Clifford biconvolution is a braided Hopf gebra 
in the usual sense? If such a braid exits how much freedom remains
for choices? Compare with [Oziewicz 1997, p. 1272] where it was shown 
that for $\dim_\R V=1$, $\dim_\R \Cl=2$ exists a 12-parameter family of 
crossings.\end{prob}

\begin{Lemma}Let $A\in\End_\R V$ and $\dim_\R V=2.$ Then the following
equations are equivalent:
\begin{align*}(i)&\quad\det(\id_V-A)\equiv\det A-\tr A+1=0,\\
(ii)&\quad(\id_V-A)\circ(\id_V+A-\tr A)=0.
\end{align*}\end{Lemma}
\begin{proof} $A^2=(\tr\,A)A-(\det A)\id_V.$ \end{proof} 
According to Lemma 5.2 we have to ask that
\begin{gather*}
\begin{array}{r@{\quad\quad}l}
\text{either} & \im(\id-A)\subset\ker(\id+A-\tr A)\\
\text{or} &\im(\id+A-\tr A)\subset\ker(\id-A).
\end{array}
\end{gather*}
We present three examples of antipode-less Clifford biconvolutions
$\Cl(\eta,\xi)$, $\det(\id-\xi\circ\eta)=0$ for $\dim_\R V=2$, $\dim_\R
V^\wedge=4$ and for $\eta\,\&\,\xi$ given by \eqref{eta} with signature
$(+,-),$
$$\
tr(\mu S)=\half\tr\{(\eta-\eta^T)\circ(\xi-\xi^T)\}=4tz,
$$
\begin{align*}
\text{Case I.}&\quad \eta^T=\eta,\quad \xi^T=\xi\;\text{and}\;
\det(\eta\circ\xi)+1=\tr(\eta\circ\xi).\\
\text{Case II.}&\quad r^2=0,\quad\tr(\xi\circ\eta)+1+tz=0,\quad
z+t\det\xi=0.\\
\text{Case III.}&\quad v^2=0,\quad tz=-1,\quad u^2=-s^2z^2.\end{align*}

\section{Splitting idempotent} 
\begin{dfn}[Eilenberg 1948, Cartan \& Eilenberg 1956] Let $R$ be a
commutative ring. An exact sequence of homomorphisms of $R$-modules,
$\im\,s=\ker\,r,$
\begin{gather}\label{sequence}
\begin{CD}0@>>>X^\prime@>{s}>>X@>{r}>>{X^{\prime\prime}}@>>>0\end{CD}
\end{gather}
splits if $\im\,s=\ker\,r=X^\prime$ is a direct summand of $X.$\end{dfn}
\begin{thm}[Cartan \& Eilenberg 1956, Scheja et al. 1980] The following
statements are equivalent
\begin{gather*}\begin{array}{r@{\quad}l}
i)  & \text{The sequence (\ref{sequence}) splits}. \\
ii) & \exists g^{\prime\prime} \in \Hom(X^{\prime\prime},X)
\quad\text{with}\;r\circ g^{\prime\prime}=\id_{X^{\prime\prime}}.\\
iii)&\exists g\in\Hom(X,X^{\prime})\;\text{with}\; g\circ s=
      \id_{X^{\prime}}.\end{array} \end{gather*}\end{thm}

\begin{dfn}[Pierce 1982; Hahn 1994]Let
$\Delta^\xi\in\alg(\Cl,\Cl\otimes\Cl)$ split the exact sequence of
$\R$-algebra homomorphisms
$$\begin{CD}0@>>>\ker\wedge^\eta@>{\Delta^\xi}>>\Cl\otimes\Cl
@>{\wedge^\eta}>>\Cl@>>>0.\end{CD}$$
Thus $\wedge^\eta\circ\Delta^\xi=\id_{\Cl}.$ In this case the element 
of the crossed $\R$-algebra, viz. $\Delta^\xi 1=\xi^\wedge\in\Cl\otimes
\Cl$, is said to be a \textit{splitting idempotent} (a cleft of
$\Cl\otimes\Cl$), $\xi^\wedge=\Delta^\xi 1=\Delta^\xi(1\cdot 1)
=(\xi^\wedge)^2$.\end{dfn}

If $\wedge^\eta\in\alg(\Cl\otimes\Cl,\Cl)$ then 
$(\Cl\otimes\Cl)\cdot(1\otimes 1-\xi^\wedge)\subset\ker\wedge^\eta.$

A crossing defined in Fig. \ref{crossing} gives a cogebra map
$\wedge^\eta\in\cog(\Cl\otimes\Cl,\Cl)$ and an algebra homomorphism
$\Delta^\xi\in\alg(\Cl,\Cl \otimes \Cl)$ [Oziewicz 1997, 2001, Cruz \&
Oziewicz 2000]. However an algebra homomorphism $\Delta^\xi$ in general
does not need to split.
\begin{dfn} A Clifford convolution $\Cl(\eta=-\eta^T,\xi=-\xi^T)$ is
said to be a Gra{\ss}mann convolution.\end{dfn} 
\begin{thm} If $\Cl(\eta,\xi)$ is a Gra{\ss}mann Hopf gebra,
$\dim_\R\Cl=4,$ then $\Delta^\xi$ splits,
$\wedge^\eta\circ\Delta^\xi=\id_{\Cl},$ iff
$$
\begin{array}{r@{\quad}l}
\text{either} & \tr(\xi\circ\eta)=\dim_\R V^\wedge,\\
\text{or}     & \eta=0\;\text{or}\;\xi=0,
\quad\text{and thus}\quad\det(\id-\xi\circ\eta)=1.
\end{array}$$\end{thm}

\section{Crossing} 
\begin{figure}[ht]
\hfil \includegraphics[height=2truecm]{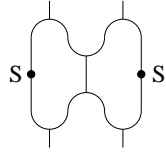}\hfil 
\caption{Definition of the crossing $\sigma$ for antipodal convolution}
\label{crossing}\end{figure}

A crossing for an antipodal convolution 
is defined on Fig. \ref{crossing}. 
The crossing $\sigma$ Fig. \ref{crossing} is equivalent that there is a
cogebra map $\wedge^\eta\in\cog(\Cl\otimes_\sigma\Cl,\Cl)$ and an
algebra homomorphism $\Delta^\xi\in\alg(\Cl,\Cl\otimes_\sigma\Cl)$
[Oziewicz 1997, 2001, Cruz \& Oziewicz 2000].

A crossing for $\Cl(0,0),$ and thus for $\tr(\mu S)=0,$ is the
involutive graded switch [Sweedler 1969, Ch.XII],
$$s(x\otimes y)\equiv(-1)^{(\text{grade}\,x)(\text{grade}\,y)}y\otimes
x,\quad s^2=\id_{\Cl\otimes\Cl}.$$

In the sequel $\mu\equiv\det(\id-\xi\circ\eta)\neq 0.$ The degree of the
minimal polynomial of the crossing $\sigma\in\End\R(\Cl\otimes\Cl)$ we
denote by: degree$(\sigma)\,\in\N.$

\begin{thm}[Oziewicz 1997, p. 1271] Let $\dim_\R V=1.$ Then 
\begin{gather*}
\det\sigma=\left(1-\frac{2}{\mu}\right)^2,\qquad\tr\,\sigma=1-\mu,\\
\begin{align*}(i)&\quad\text{degree}(\sigma)=2\quad\text{iff}\quad
\mu=1,\,2\;\text{or}\;4.\\
(ii)&\quad\text{degree}(\sigma)=3\quad\text{iff}\quad
\mu^3-4\mu^2+\mu-2=0.\\
(iii)&\quad\text{degree}(\sigma)=4\quad\text{otherwise.}
\end{align*}\end{gather*}\end{thm}

Let $\Z$ be the ring of integers, $\Z[x]$ be the ring of polynomials
in $x$ with coefficients in $\Z,$ and $p(x)\in\Z[x]$ be the following
polynomial,
\begin{eqnarray*}p(x)=&&2^3-3^2(17)x+(3)(11)x^2+(3)(107)x^3+(2)(97)x^4\\
 & &-2^2(3)(17)x^5-(5)(41)x^6-(109)x^7-(167)x^8\\
 & &+2^4x^9+2^2(3)(7)x^{10}+(37)x^{11}-(2)(3)x^{12}\end{eqnarray*} 
\begin{thm}[Degree of minimal polynomial] Let $\Cl(\eta,\xi)$ be a
Gra{\ss}mann Hopf gebra, \ie an antipodal Gra{\ss}mann convolution
(Definitions 2.3 \& 6.4) with $\tr(\mu S)=4(\sqrt{\mu}-1)$ and $\dim_\R
V=2$. The minimal polynomial of the crossing is of order 30 in 
$tr(\eta\circ\xi).$ If $tr(\eta\circ\xi)=2$ then the minimal polynomial
of $\sigma$ vanishes, as also no antipode exists in this case.
Moreover
\begin{align*}(i)&\quad\text{If}\;0\neq(\sqrt{\mu}-1)
\;\text{is not a root of a principal ideal}\;
(p(x))<\Z[x],\\
&p(\sqrt{\mu}-1)\neq 0,\;
\text{then the minimal polynomial of the crossing}\\
&\text{is of degree 3.}\\
(ii)&\quad\text{If}\;p(\sqrt{\mu}-1)=0,\;
\text{then the the minimal polynomial of the crossing is}\\
&\text{of degree}\;\leq 2.\end{align*}\end{thm}

\section*{Acknowledgement} 
CLIFFORD, a Maple V package by Rafa{\l} Ab{\l}amowicz [Ab{\l}amowicz, 
1995-2000] and extensions like BIGEBRA, by Rafa{\l} Ab{\l}amowicz and 
Bertfried Fauser, have been used to check and get some of the results.
\end{document}